 \documentclass[12pt]{amsart}
\usepackage{amssymb}
\newtheorem{theorem}{Theorem}
\newtheorem{lemma}[theorem]{Lemma}
\newtheorem{proposition}{Proposition}
\newtheorem{corollary}[theorem]{Corollary}

\theoremstyle{remark}
\newtheorem{remark}{Remark}

\newcommand{\calR}{\mathcal{R}}

\newcommand{\calD}{\mathcal{D}}

\newcommand{\bbP}{\mathbb{P}}
\newcommand{\bbQ}{\mathbb{Q}}
\newcommand{\bbR}{\mathbb{R}}
\newcommand{\bbG}{\mathbb{G}}
\newcommand{\bbZ}{\mathbb{Z}}
\newcommand{\bbC}{\mathbb{C}}

\newcommand{\bbF}{\mathbb{F}}

\newcommand{\bfF}{\mathbf{F}}

\newcommand{\la}{\langle}
\newcommand{\ra}{\rangle}

\newcommand{\Gal}{\textup{Gal}}
\newcommand{\cha}{\textup{char}}

\newcommand{\GL}{\textup{GL}}
\newcommand{\Spec}{\textup{Spec}}
\newcommand{\Pic}{\textup{Pic}}

\newcommand{\PGL}{\textup{PGL}}

\newcommand{\Hom}{\textup{Hom}}

\newcommand{\NE}{\textup{NE}}
\newcommand{\Cr}{\textup{Cr}}
\newcommand{\Aut}{\textup{Aut}}
\newcommand{\rank}{\textup{rank}}
\newcommand{\Sk}{\textup{Sk}}

\newcommand{\mori}{\overline{\textup{NE}}}

\begin{document}
 \title[On elements of prime order]{On  elements of prime order in the plane Cremona group  over a perfect field}

\author{Igor V. Dolgachev}

\address{Department of Mathematics, University of Michigan, 525 E. University Av., Ann Arbor, Mi, 49109, USA}
\email{idolga@umich.edu}
\thanks{The  first author was supported in part by NSF grant 0245203.}

\author{Vasily A. Iskovskikh}

\address{Steklov Institute of Mathematics, Gubkina 8, 119991, Moscow GSP-1, Russia}
\email{iskovsk@mi.ras.ru}
\thanks{The  second author was supported in part by RFBR 08-01-00395-a
RFBR 06-01-72017-MNTI-a, grant CRDF RUMI 2692-MO-05 and grant of
NSh-1987.2008.1}
\begin{abstract}
We show that the plane Cremona group over a perfect field $k$ of characteristic $p \ge 0$ contains an element of prime order $\ell\ge 7$ not equal to $p$ if and only if there exists a 2-dimensional  algebraic torus $T$ over $k$ such that $T(k)$ contains an element of order $\ell$. If $p = 0$ and $k$ does not contain a primitive $\ell$-th root of unity, we show that there are no elements of prime order $\ell > 7$ in $\Cr_2(k)$  and all elements of order $7$ are conjugate.
\end{abstract}

\maketitle

\bigskip\noindent
\section{Introduction}
The classification of conjugacy classes of elements of prime order
$\ell$ in the plane Cremona group $\Cr_2(k)$ over an algebraically
closed field $k$ of characteristic 0 has been known for more than a
century. The possible orders of elements not conjugate to a
projective transformation are $2,3,$ and 5 (see \cite{DI} and
historic references there). Much less is known in the case when the
field $k$ is not algebraically closed. For example, when $k = \bbQ$,
there are no elements of prime order $\ell > 3$ which are conjugate to
projective automorphisms of the plane. However, there exists an element of order
$5$ in $\Cr_2(\bbQ)$ which acts  biregulary on   a rational Del
Pezzo surface of degree 5. The first example of a birational
automorphism of $\bbP_{\bbQ}^2$ of order 7 was constructed by J.-P.
Serre \cite{Serre2}. It is realized as an automorphism of a
rational Del Pezzo surface of degree 6. He also raised the question
about the existence of birational automorphisms of prime order $\ell >
7$. As a byproduct of results of this paper, we  prove that $\Cr_2(\bbQ)$
does not contain elements of prime order $\ell > 7$  and all elements of order 7 are conjugate in $\Cr_2(\bbQ)$. 

The main result of the paper is  the following theorem.

\begin{theorem} Let $k$ be a perfect field  of characteristic $p \ge 0$. Then $\Cr_2(k)$ contains an element of prime order $\ell > 5$ not equal to $p$ if and only if there exists a 2-dimensional  algebraic $k$-torus $T$ such that $T(k)$ contains an element of order $\ell$.
\end{theorem}

Recall that an algebraic $k$-torus is an affine algebraic group $T$ over $k$ that becomes isomorphic to the torus $\bbG_m^n$ over an algebraic closure $\bar{k}$ of $k$.  In dimension 2 all $k$-tori are rational over $k$ \cite{Vos}, hence  $T(k)$  is contained in $\Cr_2(k)$. The following result of J.-P. Serre from \cite{Serre}, Theorem 4,  gives a necessary and sufficient condition in order that an algebraic $k$-torus contains an element of order $\ell$ coprime to $p$.

\begin{theorem} Let $k$ be a field of characteristic $p \ge 0$ and $\ell\ne p$ a prime  number. Then the following assertions are equivalent:
\begin{itemize}
\item[(i)] there exists a 2-dimensional algebraic $k$-torus $T$ such that $T(k)$ contains an element of order $\ell$;
\item[(ii)] the degree of the cyclotomic extension $k(\zeta_\ell)$ over $k$ takes values in  the set $\{1,2,3,4,6\}$. Here $\zeta_\ell$ is a generator of the group of $\ell$-th roots of unity in $\bar{k}$.
\end{itemize}
\end{theorem} 

We will also prove the following  uniqueness result.

\begin{theorem} Assume that $k$ is of characteristic 0 and does not contain a primitive $\ell$-th root of unity. Then $\Cr_2(k)$ does not contain elements of prime order $\ell > 7$ and  all elements of order 7 in $\Cr_2(k)$ are conjugate to an automorphism of a Del Pezzo surface of degree 6.
\end{theorem}

We are using the same tools applied in \cite{DI} for classification
of conjugacy classes of finite subgroups of $\Cr_2(\bbC)$ adjusted
to the case when the ground field is not algebraically closed. As in \cite{DI} we show that there is a bijection between the conjugacy classes of finite subgroups $G$ of $\Cr_2(k)$ and $G$-equivariant birational isomorphism classes of projective $k$-rational smooth surfaces $X$ equipped with a biregular action of $G$. 
Then we use   an equivariant version of 
two-dimensional Mori theory of minimal models over $k$. It allows us to describe minimal pairs $(X,G)$.

The authors believe that, granted time and
patience, the methods of \cite{DI} and \cite{Serre} can be used  for
classification of  the conjugacy classes of all finite subgroups of
$\Cr_2(k)$.

We thank J.-P. Serre for letting us know about his results published in \cite{Serre} that were crucial for most of results in this paper.  His numerous correspondence with us on the subject of this work  is greatly appreciated. 

\section{Minimal rational $G$-surfaces}
Let  $k$ be a  field and $\bar{k}$ be its algebraic closure. A
\emph{geometrically rational surface} $X$ is a smooth projective
surface over $k$ such that $\bar{X}$ is
birationally isomorphic to $\bbP_{\bar{k}}^2$.  Here and below, for any scheme $Z$ over $k$ we denote by $\bar{Z}$ the scheme $Z\times_k\Spec~\bar{k}$. A geometrically rational surface is called \emph{$k$-rational} if it is $k$-birational to $\bbP_k^2$.

The following
result is a culmination of several results due to V. Iskovskikh and Yu. Manin. Its modern  proof based on the theory of elementary links can be found in  \cite{Isk2}, \S 4, p.642.

\begin{theorem}\label{th1} A minimal geometrically rational surface $X$ over a perfect field $k$ is $k$-rational if and only if the following two conditions are satisfied:
\begin{itemize}
\item[(i)] $X(k) \ne \emptyset$;
\item[(ii)] $d = K_X^2 \ge 5.$
\end{itemize}
\end{theorem}

Let $G$ be a finite
subgroup of $\Aut_k(X)$. For any extension $k'/k$ the group $G$ acts
naturally on $X_{k'} = X\times_k \Spec~k'$ and hence can be considered as a
subgroup of $\Aut_{k'}(X_{k'})$. 

A \emph{geometrically rational $G$-surface} is a pair $(X,G)$, where
$X$ is a geometrically rational surface over $k$ and  $G$ is a
finite subgroup of  $\Aut_k(X)$.  A geometrically rational
$G$-surface $(X,G)$ is called minimal if any $G$-equivariant
birational $k$-morphism $X\to X'$ to a geometrically rational
$G$-surface $X'$ is a   $k$-isomorphism. In the case when $G =
\{1\}$, a minimal surface is just a $k$-minimal surface in the sense
of the theory of minimal models. If $G\ne \{1\}$, a minimal rational
$G$-surface  is not necessarily a $k$-minimal surface.

 For lack of reference,  we will prove the following theorem that is essentially due to S. Mori.

\begin{theorem}\label{thm1} Let $(X,G)$ be a minimal geometrically rational $G$-surface over a perfect field $k$. Then one of the following two cases occurs:
\begin{itemize}
\item [(i)] $X$ is a Del Pezzo surface with $\Pic(X)^G \cong \bbZ$ of degree $d:=K_X^2\neq 7$;
\item [(ii)] $X$ admits  a conic  bundle structure  with $\Pic(X)^G \cong \bbZ^2$ generated, after tensoring with $\bbQ$,  by $K_X$ and the class of a fibre.
\end{itemize}
\end{theorem}

Recall that a conic bundle structure on $X$ is a morphism $\phi:X\to
C$, where $C$ is a smooth genus 0 curve (not necessary isomorphic to
$\bbP_k^1$) and all smooth fibres over closed points $c$ are isomorphic to a reduced conic over
the residue field $k(c)$ at $s$. A non-smooth fibre over a closed point $t$ becomes isomorphic to a
bouquet of two $\bbP^1$'s over a quadratic extension of the residue
field $k(t)$ of $t$. If  $c$ denotes the number of singular geometric fibres of
$\phi$, then $d:= K_X^2 = 8-c$.

Let us first introduce the main attributes of Mori theory. These
are the space $N(X) = \Pic(X)\otimes \bbR$,  the cone
$\overline{NE}(X)\subset N(X)$ spanned by  effective divisor
classes in $\Pic(X)$ and the part $\overline{NE}(X)_{K_X\ge 0}$  of
$\overline{NE}(X)$ defined by the inequality $K_X\cdot D \ge 0$.

Let $N(\bar{X})$, $\overline{NE}(\bar{X})$ and
$\overline{NE}(\bar{X})_{K_{\bar{X}}\ge 0}$ be the similar objects
for the surface $\bar{X}$. The Galois group  $\Gamma_k =
\text{Gal}(\bar{k}/k)$ of $k$ acts naturally on $\bar{X}$ and hence
acts on $N(\bar{X})$ and $\overline{NE}(\bar{X})$. The action
factors through the action of a finite quotient of $\Gamma_k$ and we
have (\cite{Mori}, Proposition 2.6).
$$N(\bar{X})^{\Gamma_k} = N(X), \ \overline{NE}(\bar{X})^{\Gamma_k} = \overline{NE}(X), \ \overline{NE}(\bar{X})_{K_{\bar{X}}\ge 0}^{\Gamma_k} = \overline{NE}(X)_{K_{X}\ge 0}, $$
\begin{equation}\label{cone1}
\mori(X) = \mori(X)_{K_X\ge 0}+\bbR_{\ge 0}[\ell_1]+\cdots+\bbR_{\ge 0}[\ell_s],
\end{equation}
where $\ell_1,\ldots,\ell_s$ are $k$-irreducible curves such that
$\bar{\ell}_i $ is equal to the
$\Gamma_k$-orbit of  an extremal curve, i.e.  a smooth rational
curve $\bar{E}_i$ such that $0 >K_{\bar{X}}\cdot \bar{E}_i \ge
-3$.

By Theorem 2.7 from \cite{Mori} one of the following cases occurs
\begin{itemize}
\item[(i)] there exists $\bar{\ell}_i$ equal to  a disjoint sum of $(-1)$-curves on $\bar{X}$. In this case   there exists a birational $k$-morphism $f:X\to Y$ onto a smooth surface $Y$ such that $f(\ell_i)$ is a point.
\item[(ii)] $\text{rank}~\Pic(X) = 2$, $\ell_1^2 = 0$ and there exists a conic bundle structure  on $X$ with the class of a fibre over a closed point equal to the divisor class of $\ell_1$.
 \item[(iii)] $\Pic(X) \cong \bbZ$ and $X$ is a $k$-minimal Del Pezzo surface of degree $ d\ne 7$.
 \end{itemize}
Applying the averaging operator $\sum_{g\in G}g$ to  \eqref{cone1},
we obtain
\begin{equation}\label{cone2}
\mori(X)^G = \mori(X)_{K_X\ge 0}^G+\bbR_{\ge 0}\big[\sum_{g\in
G}g(\ell_1)\big]+\cdots+\bbR_{\ge 0}\big[\sum_{g\in G}(g(\ell_s)\big].
\end{equation}
Write $ \bbR_{\ge 0}\big[\sum_{g\in G}g(\ell_i)\big] = \bbR_{\ge 0}[L_i]$,  where $L_i$
is a reduced $G$-orbit of $\ell_i$. Let $L$ be one of the orbits which  generates  an extremal ray in
$\mori(X)^G$ (since $K_X$ is not nef, it exists).  

From now on we assume that $\text{rank}~\Pic(X)^G
\ge 2$, otherwise the assertion  of Theorem \ref{thm1} is obvious. Note that the case $d = 7$ is excluded for the obvious reason that $\bar{X}$ contains three $(-1)$-curves forming a chain, the middle one is $\Gamma_k\times G$-invariant.

The proof of  Lemma 2.5 from \cite{Mori}, shows that $L^2 \le 0$
since otherwise $[L]$ lies in the interior of  $\mori(X)^G$.  Since $\bbR_{\ge 0}[L]$ is an extremal
ray in $\mori(X)^G$, the divisor $\bar{L}$ is a $\Gamma_k\times
G$-orbit of an irreducible component  $\bar{E}$ of $\bar{\ell}$.

Assume $\bar{E}^2 = 0$ as in case (ii) from above. Then $\bar{L}^2
\le 0$ implies that $\bar{L}$ is the disjoint sum of irreducible
smooth rational curves with self-intersection $0$. By the Hodge
index theorem, they are proportional in $\mori(\bar{X})$ (otherwise the signature of the intersection form on $\Pic(\bar{X})$ is not equal to $(1,r)$) and the
contraction of $\Gamma_k\times G$-invariant extremal ray
$\bbR_{\ge 0}[\bar{L}]$
defines a $G$-equivariant $k$-map to some $k$-curve $C$ with
geometrically connected fibres (see \cite{Mori}, 2.5.1).

Assume we are in the case (i) from above, i.e. $\bar{E}$ is a $(-1)$-curve. Let   $\bar{L}= \bar{E}_1+\cdots+\bar{E}_m,$
where $\{\bar{E}_1,\ldots,\bar{E}_m\}$  is the orbit of $ \bar{E} = \bar{E}_1$ with respect to  $\Gamma_k\times G$. We have
\begin{equation}\label{zhang} 0\ge \bar{L}^2 = m\bar{E}_1\cdot (\bar{E}_1+\cdots+\bar{E}_m) = m(-1+\sum_{i=2}^m\bar{E}_1\cdot \bar{E}_i).\end{equation}
If $\bar{L}^2 < 0$, then $\sum_{i=2}^m\bar{E}_1\cdot \bar{E}_i = 0$.
Replacing $\bar{E}_1$ with some other $\bar{E}_i$, we obtain that
$\bar{L}$ is a disjoint sum of $(-1)$-curves.
Thus $L$ can be $G$-equivariantly contracted contradicting the
minimality of $(X,G)$.

So we may assume that $\bar{L}^2 = 0$ and
$\sum_{i=2}^m\bar{E}_1\cdot \bar{E}_i = 1$. Without loss of
generality, $\bar{E}\cdot \bar{E}_2 = 1, \bar{E}\cdot \bar{E}_i = 0,
i > 2.$ Write  $\bar{E}_2 = g(\bar{E})$ for some $g\in
\Gamma_k\times G$. Then $\bar{E}\cdot \bar{E}_2 =
g^{-1}(\bar{E})\cdot \bar{E} = 1$, hence $g^{-1}(\bar{E}) =
\bar{E}_2 = g(\bar{E})$ and $g^2$ leaves $\bar{E}$ invariant. If $m
= 2$, we obtain that the linear system $|\bar{E}_1+\bar{E}_2|$
defines a $G$-invariant $k$-map onto $\bbP^1$. Assume $m \ge 3$.
Replacing $\bar{E}_1$ with $\bar{E}_3$ in \eqref{zhang}, and
repeating  the argument, we may assume that $\ell$ is equal to the
disjoint sum of $\frac{m}{2}$-bouquets $B_i$ of two $(-1)$-curves, and each $g\in \Gamma_k\times G$ either leaves a
bouquet invariant or sends it to another bouquet. Since $B_i^2 = 0,
B_i\cdot K_X = -2$, and $B_i\cdot B_j = 0, i\ne j$, by the Hodge index
theorem, the divisors $B_i$ represent the same divisor class. The
linear system $|nB_i|$, for some $n > 0$, defines a $G$-invariant
$k$-map to a curve $C$ with geometrically connected fibres.

It remains to show that a $G$-invariant $k$-map $f:X\to C$
constructed as above is a conic bundle structure and also
$\text{rank}~\Pic(X)^G = 2$. The divisor class $L$ defines a $G$-equivariant extremal ray. Thus its contraction from above has the relative Picard rank equal to 1 defining a relatively minimal conic bundle.  This implies that  $\text{rank}~\Pic(X)^G = 2$. 

\begin{remark} The assumption that the ground field $k$ is perfect is not essential. We used it only to simplify the notations (see \cite{Mori}, Theorem 2.7, where $k$ is not assumed to be perfect).
\end{remark}

The following lemma is well-known. Its proof in \cite{DI}, Lemma 3.5 can be extended to the case of an arbitrary  ground field $k$ by referring to \cite{Lipman} where one can find a proof of  the existence of a $G$-equivariant resolution of surfaces. In fact, in \S  2 of his paper Lipman proves that the known existence of any resolution of a normal surface $X$ implies that the sequence
$$X = X_1 \leftarrow X_2  \leftarrow \ldots  \leftarrow X_n  \leftarrow \ldots ,$$
where each morphism $X_{i+1}\to X_i$ is the composition of the
blow-up of all singular points and the normalization, terminates in a
nonsingular surface $Y$. Obviously, the composition map $Y\to X$ is
$G$-equivariant for any finite group $G$ acting on $X$.

\begin{lemma} Let $G$ be a finite subgroup of $\Cr_2(k)$, then there exists a $k$-rational nonsingular surface $X$,   an injective  homomorphism $\rho:G\to \Aut_k(X)$ and a birational $G$-equivariant $k$-map  $\phi:X\to \bbP_k^2$ such that
$$G = \phi\circ \rho(G)\circ \phi^{-1}.$$
\end{lemma}

The $G$-surface $(X,\rho(G))$ can be obviously replaced by a minimal
rational $G$-surface, and in this way we establish a natural
bijective correspondence between the conjugacy classes of finite
subgroups $G$ of $\Cr_2(k)$ and birational isomorphism classes of minimal nonsingular $k$-rational
$G$-surfaces $(X,G)$ (see Theorem 3.6 in \cite {DI}). 

In this note we will be interested in the  case  $k$ is any perfect field
and surfaces $X$  are $k$-rational.  Since $k$ is perfect, the surface $\bar{X}$ is a nonsingular surface over $\bar{k}$, so that we can apply the theory of nonsingular rational surfaces over an algebraically closed field.

 It follows from above that we may assume that an automorphism $\sigma$ of prime order $\ell$ acts on a $k$-rational surface $X$ making the pair $(X,\la \sigma\ra)$ a minimal $\la \sigma\ra$-surface. We say that $\sigma$ \emph{acts minimally} on $X$.

\section{Elements of finite order in reductive algebraic groups}\label{sect3}

Here we give a short exposition of some of Serre's results from \cite{Serre} that are relevant to our discussion.

For any integer $N$ and a prime number $\ell$ we denote by $\nu_\ell(N)$ the largest $n$ such that $\ell^n$ divides $N$. For any finite group $A$ we set $\nu_\ell(A)$ to be equal to $\nu_\ell(|A|)$. 

The first historical result on elements of finite order in linear groups is the following theorem of Minkowski  \cite{Mink}.

\begin{theorem}\label{mink} Let $n$ be an integer $\ge 1$ and $\ell$ be a prime number. Define
$$M(n,\ell) = \Bigl[\frac{n}{\ell-1}\Bigr]+\Bigl[\frac{n}{\ell(\ell-1)}\Bigr]+\Bigl[\frac{n}{\ell^2(\ell-1)}\Bigr]+\cdots$$
Let $G$
 be a finite subgroup of $\GL_n(\bbQ)$. Then  $\nu_\ell(G)\le M(n,\ell)$ and there exists an $\ell$-subgroup $G$ of $\GL_n(\bbQ)$ with 
 $\nu_{\ell}(G) = M(n,\ell)$.
\end{theorem}

We will be interested in a similar result for the group $\PGL_n(k) = \GL_n(k)/k^*$. We cannot apply the previous theorem since an element of finite order does not necessary lift  to an element of finite order in $\GL_n(k)$. However we can apply the following result of Serre from \cite{Serre}.

\emph{From now on $\ell$ will denote a prime odd number not equal to the characteristic of $k$}.

Let $\zeta_\ell$ be a generator of the group $\mu_\ell(\bar{k})$ of $\ell$-th roots of unity in $\bar{k}$. Set 
$$t_\ell = [k(\zeta_\ell):k],\quad m_\ell = \sup\{d\ge 1:\zeta_{\ell^d}\in k(\zeta_\ell)\}.$$
Let $P$ be the prime field contained in $k$. By Galois Theory,
$$t_\ell = [P(\zeta_\ell):P(\zeta_\ell)\cap k].$$
If $\cha~k = 0$, we obtain that  $t_\ell$ divides $\ell-1$. If  $k$ is of characteristic $p > 0$, then 
$P(\zeta_\ell) \cong \bbF_{p^s}$, where $s$ is the order of $p$ in $\bbF_\ell^*$. Thus $t_\ell$ divides $s$.

For example, when $k = \bbQ$, we have $t_\ell = \ell-1$ and $m_\ell = 1$. If $k = \bbF_{q}$, then $t_\ell $ is equal to the order of $q$ in $\bbF_\ell^*$ and $m_\ell = \nu_\ell(q^{\ell-1}-1)$.

The following  is a special case of Theorem 6 from \cite{Serre}.

\begin{theorem}\label{serre1} Let $A$ be a finite subgroup of  $\PGL_{n+1}(k)$. For any $\ell > 2$, 
$$\nu_\ell(A) \le \sum_{2\le s\le n+1, t_\ell|s}(m_\ell+\nu_\ell(s)).$$
\end{theorem}

\begin{corollary}\label{no2} Assume $t_\ell \ge n+2$. Then $\PGL_{n+1}(k)$ does not contain elements of prime order $\ell$.
\end{corollary}

For example, if   $k$ is of characteristic zero and $m_\ell(k) =\{1\}$ (e.g. $k = \bbQ$), then $t_\ell = \ell-1$ and we get 
$$n \ge \ell-2$$
 if $\PGL_{n+1}(k)$ contains an element of order $\ell$.  For example, $\PGL_{n+1}(k)$ contains an element of order 7 only if $n\ge 5$.

On the other hand, if $k = \bbF_2$ and $\ell = 7$,  then $t_\ell = 3$ and  it is known that $\PGL_3(k)$ is isomorphic to a simple group of order 168 and it contains an element of order 7. 

The next result of Serre \cite{Serre}, Theorems 4 and 4', concerns elements of finite order in an algebraic $k$-torus. 

\begin{theorem}\label{serre2} Let $T$ be an algebraic $k$-torus and $A$ be a finite subgroup of $T(k)$. Then 
$$\nu_\ell(A) \le m_\ell\Big[ \frac{\dim T}{\phi(t_\ell)}\Big],$$
where $\phi$ is the Euler function.
Assume $m_\ell < \infty$ (e.g. $k$ is finitely generated over its prime subfield).  For any $n\ge 1$ there exists an $n$-dimensional $k$-torus $T$ and a finite subgroup $A$ of $T(k)$ such that $\nu_\ell(A) = m_\ell\big[ \frac{\dim T}{\phi(t_\ell)}\big].$
\end{theorem}

Since this will be of importance to us, let us give Serre's construction proving the last assertion. 

Let $T$ be an algebraic $k$-torus  of dimension $d$ over an arbitrary  field $k$. By definition, it is an affine algebraic group over $k$ that becomes isomorphic to the group $\bbG_{m}^d $ over some finite extension $E$ of $k$.  One can always choose $E$ to be a separable Galois extension, it is called  a \emph{splitting field} of $T$. 

We denote by $M = \Hom(\bbG_{m}^d,\bbG_m)$ the group of rational characters of  $\bbG_{m}^d$, so that  $\bbG_{m,E}^d \cong \Spec~\bbZ[M]$. A $d$-dimensional $k$-torus $T$ split over a Galois extension $E/k$ with Galois group $\Gamma$ defines a structure of a  $\Gal(E/k)$-module on $M$ such that $T \cong \Spec~E[M]^{\Gamma}$. In this way the category of $d$-dimensional $k$-tori split over $E$ is anti-equivalent to the category of free  abelian groups equipped with a structure of a $\Gamma$-module \cite{Vos}. 

The group algebra  $\bbZ[\Gamma]$  equipped with the natural structure of a $\bbZ[\Gamma]$-module $M$ (the regular representation of $\Gamma$) defines an  algebraic $k$-torus $R_{E/k}(\bbG_m)$ of dimension equal to $[E:k]$. It represents the functor on the category of $k$-algebras defined by $L \to (L\otimes_kE)^*$. 

Now let us take $E = k(\zeta_\ell)$ and let $\Gamma$ be the Galois  group of $E/k$. It is a cyclic group of order  $t = t_\ell$. Let 
$$x^{t} -1 = \Phi_{t}(x)\Psi(x),$$
where $\Phi_{t}(x)$ is the $t$-th cyclotomic polynomial over $\bbQ$. We have $\bbZ[\Gamma] \cong \bbZ[x]/(x^t-1)$. The multiplication by $\Psi(x)$ in $\bbZ[x]$ defines an inclusion of $\Gamma$-modules $\bbZ[x]/(\Phi_t(x)) \hookrightarrow \bbZ[x]/(x^t-1)$ and hence  a surjective homomorphism of $k$-tori $R_{E/k}(\bbG_m)\to T$. Here  $T$ is a $\phi(t)$-dimensional $k$-torus defined by the character module $M = \bbZ[x]/\Phi_t(x)$.  One can view $T$ as the image of $R_{E/k}(\bbG_m)$ under the endomorphism $\Psi(\gamma)$, where $\gamma$ is a generator of $\Gamma$. One can show that $\Psi(\gamma)$ acts  on the subgroup $\la \zeta_{\ell^{m_\ell}}\ra$ of $E^*$ by an automorphism. This implies that  $T$ contains  $\zeta_{\ell^{m_\ell}}$.  Thus 
 $$ m_\ell\Big[ \frac{\dim T}{\phi(t)}\Big] =  m_\ell\Big[\frac{\phi(t)}{\phi(t)}\Big] = m_\ell = \nu_\ell(\la\zeta_{\ell^{m_\ell}} \ra).$$

Now, for any positive integer $n$ let $s = \big[\frac{n}{\phi(t)}\big]$. Take $T$ equal to $T_1^s\times \bbG_{m}^{n-s\phi(t)}$, where $T_1$ is constructed as above. It is easy to see that $T$ satisfies the  last assertion of Theorem 10.
 
\begin{corollary}\label{cor} A 2-dimensional $k$-torus $T$ with $T(k)$ containing an element of prime order $\ell > 2$ exists if and only if $t_\ell$ takes values in the set $\{1,2,3,4,6\}$.
\end{corollary}

\begin{proof} In fact, the set  $\{1,2,3,4,6\}$ is the set of positive integers $t_\ell$ such that $\phi(t_\ell)\le 2$. If $\phi(t_\ell) > 2$, Serre's bound implies that no such torus exists. If $\phi(t_\ell) = 2$, Serre's construction from above exhibits such a torus. If $\phi(t_\ell) = 1,$ i.e. $t_\ell = 1$ or $2$, we can take $T = \bbG_{m,k}^2$ in the first case and $T = R_{k(\zeta_\ell)/k}(\bbG_{m})$ in the second case.
\end{proof}

\section{Del Pezzo surfaces of degree 6}

Here we recall some well-known facts about  toric Del Pezzo surfaces of degree 6.  Let $\bar{k}$ be an algebraic closure of $k$. We continue to assume that $k$ is perfect. Recall that a Del Pezzo surface $S$ of degree 6 over $\bar{k}$ is isomorphic to the blow-up of 3 non-collinear points $p_1,p_2,p_3$ in $\bbP_{\bar{k}}^2$. The set of $(-1)$-curves on $S$ consists of six curves, the exceptional curves of the blow-up morphism, and the proper transforms of the lines $\overline{p_i,p_j}$. In the anti-canonical embedding $S\hookrightarrow \ \bbP_{\bar{k}}^6$, they are six lines forming a hexagon.

The surface $S$ is isomorphic to $\calD_{\bar{k}} = \calD\otimes \bar{k}$, where $\calD$ is a unique smooth projective toric surface (defined over $\bbZ$) with the Picard group of rank $4$. It is defined  by a complete $\bbZ^2$-fan $\Sigma$ whose $1$-skeleton $\Sk^1(\Sigma)$ consists of 6 rays, the primitive vectors in $ \bbZ^2$ spanning these rays are $\pm e_1,\pm e_2, \pm (e_1+e_2)$, where $e_1,e_2$ is a basis in $\bbZ^2$. The closures of the orbits corresponding to these rays  is the set of  six $(-1)$-curves. 

The group of  automorphism of the surface $\calD$ is a group scheme over $\bbZ$  isomorphic to $\bbG_m^2\rtimes D_{12}$, where $D_{12}$ is the dihedral group of order 12,  realized as the subgroup of $\GL_2(\bbZ)$ leaving the fan $\Sigma$ invariant.

Let $X$ be a rational surface Del Pezzo surface of degree 6 over $k$. Then   $\bar{X} = X\otimes_k\bar{k}$ is isomorphic to $\overline{\calD} = \calD\otimes\bar{k}$. Since the set of all $(-1)$-curves on $\bar{X}$ is defined over $k$, its complement $U$ in $X$ becomes isomorphic to a torus over $\bar{k}$. This implies that $U$  is a torsor (= principally homogeneous space) over a two-dimensional $k$-torus $T$ (see \cite{Manin}, Chapter IV, Theorem 8.6). Since $X$ is rational, $X(k)\ne \emptyset$ and hence $U(k)\ne \emptyset$ (\cite{KV}, Proposition 4). This shows that $U$ is an algebraic $k$-torus.

Fix a structure of $k$-torus  $T$ on $U$ and consider its action on itself by translation. Then we can extend the action to the action of $T$ on $X$ (see \cite{Rob}, p. 22). Thus $X$ is a (not necessarily split)  toric variety over $k$. Let $E$ be a splitting field of the torus $T$ and $X_E = X\otimes_kE$. Then the split torus $T_E$ acts on $X_E$ making it into a split toric variety over $E$. It is obviously isomorphic to $\calD_E := \calD\otimes_kE$. Thus $E$ is a splitting field of $X$, i.e.  $X_E\cong \calD_E$.  Conversely, if $E$ is a splitting field of $X$, the torus $T_E$ is split. The Galois group $\Gamma$ of $E/k$ acts on $X_E$ via automorphisms of  the toric variety, i.e. it acts on the lattice $N\cong  \bbZ^2$, the dual lattice of the lattice $M$ of characters of $T$,  leaving invariant the fan $\Sigma$ defining the toric variety $\calD$. The variety $X$ is isomorphic to the descent of $\calD$ defined by this action, defined uniqueley up to $k$-isomorphism (see \cite{Vos}, \cite{BT}).

Choose a splitting field  $E$ such that the map $\Gamma \to \Aut(\Sigma)$ is injective. It follows from the classification of subgroups  of  $D_{12}$ that one of the following cases occurs.

\begin{itemize}
\item[(i)] $\Gamma = \la \gamma\ra, |\Gamma| = 2$, the action of $\Gamma$ on $N$ is given by  $\gamma:(m,n)\mapsto (n,m)$;
\item[(ii)] $\Gamma = \la \gamma\ra, |\Gamma| =  2$,  $\gamma:(m,n)\mapsto (-m,-n)$;
\item[(iii)] $\Gamma = \la \gamma\ra, |\Gamma| = 3$, $\gamma:(m,n) \mapsto (-n,m-n)$;
\item[(iv)] $\Gamma = \la \gamma_1,\gamma_2\ra, \Gamma\cong (\bbZ/2)^2$,  $\gamma_1:(m,n) \mapsto (n,m), \gamma_2:(m,n)\mapsto (-m,-n)$;
\item[(v)] $\Gamma = \la \gamma_1,\gamma_2\ra, \Gamma \cong  S_3$,  $\gamma_1:(m,n) \mapsto (-n,m-n), \gamma_2:(m,n)\mapsto (n,m)$;
\item[(vi)] $\Gamma = \la \gamma\ra, \Gamma \cong \bbZ/6$, $\gamma:(m,n) \mapsto (m-n,m)$;
\item[(vii)] $\Gamma = \la \gamma_1,\gamma_2\ra, \Gamma = D_{12}$,  $\gamma_1:(m,n) \mapsto (m-n,m), \gamma_2:(m,n)\mapsto (n,m)$.
 \end{itemize}
We easily get (see \cite{BT}) that 
$$\rank~\Pic(X) =\begin{cases} 3& \text{in cases  (i), (ii)},\\
2&\text{in cases (iii),(iv), (v)},\\
1&\text{in cases (vi),(vii)}.
\end{cases}
$$
Note that in all cases except case (i), the torus $T$ is anisotropic (i.e. $M^{\Gamma} = \{0\}$). 

\begin{proposition}\label{dp6} Assume a cyclic group $G = \la \sigma\ra$ of   prime order $\ell \ge 5$ acts minimally on a $k$-rational Del Pezzo surface $X$  of degree 6. Let $T$ be the complement of the union of $(-1)$-curves on $X$ that acts on $X$ via its structure of a toric surface over $k$. Then $\sigma$ is defined via the action by an element $\tilde{\sigma}  \in T(k)$.  The torus $T$  splits  over   $k(\zeta_\ell)$ with cyclic Galois group $\la \gamma\ra$ of order 6. The $G$-surface $(X,G)$ is unique up to $k$-isomorphism. 
\end{proposition}

\begin{proof} By assumption $\ell \ge 5$, the group $G$ cannot be isomorphic to a subgroup of $D_{12}$, hence can be identified with a subgroup of $T(k)$. The assumption that $(X,G)$ is minimal implies that $X$ is $k$-minimal because the set of $(-1)$-curves does not contain orbits of length $\ell$. Let $E$ be a minimal splitting field of $X$ (and hence for $T$) with Galois group $\Gamma$.  Then $\Gamma$ is isomorphic to a subgroup of $D_{12}$ and only cases (vi) and (vii) from the above list are possible. Indeed, in the remaining cases $\rank \Pic(X) > 1$ and hence $X$ is not $k$-minimal. 

Let us consider case (vi). In this case $\Gamma = \la \gamma\ra$ is a cyclic group of order 6 and $\gamma$ acts on the lattice $N = \bbZ^2$ by a matrix $\left(\begin{smallmatrix}1&-1\\0&1\end{smallmatrix}\right)$. This matrix satisfies the cyclotomic equation $x^2-x+1 = 0$. The lattice $M$ of characters of $T$ is dual to the lattice $N$, and hence it is isomorphic as a $\Gamma$-module to $\bbZ[x]/(x^2-x+1)$. Thus the torus $T$ and the surface $X$ are determined uniquely by the extension $E/k$. Let us show that, under our assumptions, $E = k(\zeta_\ell)$ and hence $t_\ell = 6$.

Since $G\subset T(k)$ extends the translation action on $T$ by an element of order $\ell$, we have $T(E) \cong (E^*)^2$ contains $G$. The image of $G$ under the projection to $E^*$ contains the group $\mu_\ell$ of $\ell$-th roots of 1. Therefore, $E' = k(\zeta_\ell)$ is contained in $E$. Let $T(E)[\ell] \cong M\otimes M/\ell M\otimes \mu_\ell$ be the $l$-torsion subgroup of $T(E)$. The group  $\Gal(E/E')$ acts trivially on $T(E)[\ell]$, or, equivalently, the actions of $\Gamma$ and $\Gal(E'/k)$ on $T(E)[\ell]$ coincide. By  Minkowski's lemma (\cite{Serre}, Lemma 1), the natural homomorphism $\GL(M)\to \GL(M/\ell M)$ is a bijection on the set of elements of finite order. This shows that the actions of the Galois groups coincide on the whole character lattice $M$ of $T$. In view of minimality of the splitting field, we get $E = E'$.

The previous argument also shows that case (vii) does not occur. In fact, if $E$ is a minimal splitting field with Galois group isomorphic to $D_{12}$, then it follows from above that $E' = k(\zeta_\ell)$ is contained in $E$ and must coincide with $E$. Since the extension $E'/E$ is cyclic this is impossible.
 \end{proof}

Applying Corollary \ref{cor}, we obtain that $\ell = 7$ if  $\cha~k = 0$  and $k\cap \bbQ(\zeta_\ell) = \bbQ$. On the other hand, if  for example
$k = \bbF_{17} $ and $\ell = 13$, we get $ t_\ell = 6$.  This shows that  $\ell = 13$ is possible in this case.

\begin{remark}\label{toric} Another example of a non-split toric Del Pezzo surface is a Del Pezzo surface $X$ of degree 8 isomorphic to a nonsingular quadric $\bfF_0$ over $\bar{k}$. Recall that $\bfF_0$ can be given by a fan $\Sigma$ with 1-skeleton spanned by vectors $\pm e_1,\pm e_2$. The automorphism group of the toric  surface (i.e. the subgroup of $\Aut(\bfF_0)$ preserving the toric structure) is equal to $\bbG_m^2\rtimes D_{8}$, where $D_8$ is the dihedral group of order 8.  Fix a toric structure on $\bbF_0$ by fixing a quadrangle of rulings and a structure of a torus on its complement. Let $E/k$ be a Galois extension  with Galois group $\Gamma$. Chose an injective homomorphism $\Gamma\to \GL(2,\bbZ)$ that leaves invariant the fan $\Sigma$. By a theorem of Voskresenki\u{i} mentioned in above, there is a unique descent of the toric variety $\bbF_0$ to a toric $k$-variety $X$. This is a toric Del Pezzo surface of degree 8. Since all 2-dimensional $k$-tori are rational, the surface $X$ is a $k$-rational.  

Assume $t_\ell\ge 3$, following the proof of the previous proposition,  we see that an automorphism of order $\ell\ge 5$ can act minimally on $X$ only if $T$ is isotropic with $\Gamma$ isomorphic to a cyclic group of order 2 acting by $(m,n)\mapsto (n,m)$ or anisotropic with $\Gamma$ isomorphic to a cyclic group of order 4 or $D_{8}$. In the first case $T = R_{E/k}(\bbG_m)$ and $T(k) = E^*$ does not contain  an element of order $\ell$.  In the second case $T_E$ is the quotient of the torus of $R_{E/k}(\bbG_m)$ corresponding to the submodule  $\bbZ[x]/(x^2+1)$ of $\bbZ[\Gamma]\cong \bbZ[x]/(x^4-1)$. The group  $T(k)$ contains an element of order $\ell\ge 5$ only if $t_\ell = 4$. If $\Gamma \cong D_8$, we obtain that $E$ contains $k(\zeta_\ell)$. Since $k(\zeta_\ell) = E$ this case is impossible.
\end{remark}

\section{The main theorem}

Let us start proving our main result which is  Theorem 1 from Introduction. Let $\sigma$ be an element of prime order $\ell > 2$ in $\Cr_2(k)$.  As before we assume that $\ell \ne p = \cha~k$.

\begin{proposition}\label{prop1} Assume $\ell \ge 5$ and $\sigma$ acts  minimally on a $k$-rational conic bundle $X$. Then 
$t_\ell \le 2$ and   $\sigma$ is conjugate in $\Cr_2(k)$ to an element defined by a rational point  on a 2-dimensional algebraic $k$-torus.
\end{proposition}

\begin{proof} Let $\pi:X\to C$ be a conic bundle structure on $X$. Since $X$ is $k$-rational, the set $X(k)$ is not empty and hence the set $C(k)$ is not empty. Thus $C\cong \bbP_k^1$.

Assume $t_\ell > 2$, by Corollary \ref{no2}, $\Aut_k(C) \cong \PGL_2(k)$ does not contain elements of  order $\ell$. Thus $G = \la \sigma \ra$ acts trivially on the base $C$ and hence acts on the generic fibre $X_\eta$ by $k(\eta)$-automorphisms. However, $t_\ell(k(\eta)) > 2$ so we can apply Corollary \ref{no2} again and get that the action on $X_\eta$, and hence on $X$, is trivial. 
 
 Now we assume that $t_\ell\le 2$. If $(\bar{X},G)$ is  minimal, then $G$ contains an element that switches irreducible components of some reducible fibre (\cite{DI}, Lemma 5.6, the proof of this lemma works whenever the ground field is algebraically closed and its characteristic is coprime with the order of the group). Thus the cyclic group $G$ of order $\ell$ is mapped nontrivially  to a group of even order. Since $\ell $ is odd this is impossible. Since $(X,G)$ is minimal, this implies that $X\to C$ is relatively minimal over $k$. Suppose $X$ is not $k$-minimal, then there exists a birational $k$-morphism from $X$ to a rational $k$-minimal Del Pezzo surface $X'$.  In terminology of \cite{Isk2} this makes a link of type I. Applying Theorem 2.6 (i) from loc.cit., we see that the possible values of $K_X^2$ are $8, 6, 5,$ and $3$. 
 
 If $K_X^2 = 8$, we have   $X' = \bbP_k^2$ and $X \cong \bfF_1$, so $(X,G)$ is not minimal.
  
 If $K_X^2 = 6$, $X'$ is a quadric $Q$ in $\bbP_k^3$ with $\Pic(Q) \cong \bbZ$. The morphism $X\to X'$ blows down a pair $E$ of conjugate $(-1)$-curves. We have $E\in |-K_X-2F|$, where $F$ is the class of a fibre in the conic bundle. Obviously $G$ leaves $E$ invariant and hence $(X,G)$ is not minimal.

 If $K_X^2 = 5$, $X' = \bbP_k^2$, and $X\to X'$ blows down a fourtuple $E$ of $(-1)$-curves.  We have $E\in |-2K_X-3F|$ and again $(X,G)$ is not minimal.
 
 Finally, if $K_X^2 = 3$, then $X$ is a Del Pezzo surface isomorphic to a cubic surface. In this case a $(-1)$-curve is a line on $X$.  Thus $X$ contains  a line $L$ defined over $k$ such that the conic bundle is formed by the pencil of planes through the line. The line $L$ is $G$-invariant and the pair $(X,G)$ is not minimal.
 
Thus we may assume that $X$ is $k$-minimal. Applying Theorem \ref{th1}, we obtain  $K_X^2 \ge 5$. The conic bundle $\bar{X} $ contains $c= 8-K_X^2\le 3$ reducible fibres. Since $G$ in its action on $C$ has no orbits of length $2$ or $3$, $G$ has at least $c$ fixed points on $C$. 

If $c = 3$, $G$ acts identically on $C$. The surface $\bar{X}$ admits a birational morphism  to $ \bbP_{\bar{k}}^2$, the blow-up of 4 points.  The exceptional curves are  sections of the conic bundle.  Their union $E$  belongs to the linear system 
$|-2K_{\bar{X}}-3F|$, where $F$ is the class of  a fibre of the conic bundle on $\bar{X}$. This shows that  $G$ leaves $E$ invariant, and, since the order of $G$ is an odd prime number $\ell$, it leaves each component invariant. Thus $G$ has 4 fixed points on each fibre, hence acts trivially on the general fibre, hence trivially on $X$

Assume $c = 2$, i.e. $K_X^2 = 6$ and $\bar{X}$ is obtained by blowing up two points on a quadric. They   are the base points of a pencil of conic plane sections of the quadric. The exceptional divisor $E$ is a 2-section of the conic bundle that belongs to the linear system 
$|-K_{\bar{X}}-2F|$. As in the previous case $G$ fixes each component of $E$. Let $T$ be the complement of the union $B$ of  $E, F_1,F_2$. Since $(X,G)$ is minimal, the Galois group switches the components of $E$ and the components of the fibres. This shows that $B(k)$ consists of two points, the singular points of the fibres. If $k$ is infinite, $X(k)$ is dense in $X(\bar{k})$ because $X$ is rational. This shows that $T(k)\ne \emptyset$. If $k = \bbF_q$ is finite, it follows from the Weil Theorem (\cite{Manin}, Chapter IV, \S 5, Corollary 1) that $\#X(k) \equiv 1 \mod q$, hence $\#X(k) > 2$ and $T(k)\ne \emptyset$  again. Thus $T$ is a $G$-invariant 2-dimensional torus. Since $\sigma$ is of order $\ell \ge 5$, its image in the automorphism group of $T$ (as an algebraic group) is trivial. Thus $\sigma$ is realized by an element of $T(k)$.  

Assume $c = 1$. Then $X$ is obtained by blowing up one point on a minimal ruled surface. This point must be a $k$-point and the exceptional curve is a component of a unique singular fibre. Since $\ell $ is odd it is $G$-invariant. Thus $(X,G)$ is not minimal.

Finally assume that $c = 0$, i.e. $\bar{X}$ is a minimal ruled surface $\bfF_e$. Since the base $C$ of the fibration is isomorphic to $\bbP_k^1$, this easily implies that $X \cong \bfF_e$ over $k$. Let $\bar{G}$ be the image of $G$ in $\Aut_k(C)$. 

Suppose $\bar{G}$ is trivial, $G$ has 2 fixed points, maybe conjugate on the general fibre isomorphic to $\bbP_{k(C)}^1$. This defines a 2-section $S$ of the fibration fixed by $G$. If $e > 0$, one of the components of $\bar{S}$ is the exceptional section. The other component must be defined over $k$ and $G$-invariant. Since $G$ acts identically on the base, each component is fixed by $G$ pointwisely. Since the fixed locus of $G$ is smooth (because $\ell \ne \cha~k$), the two  components do not intersect.  The complement of $S$ and  two  fibres is a $G$-invariant $k$-torus $T$. If $e = 0$, then $G$ leaves invariant the other projection to $\bbP_{\bar{k}}^1$. There will be invariant fibres over two points, maybe conjugate over $k$. The complement of these two fibres and two fibres of the first projection is an algebraic torus defined over $k$. 

 If $\bar{G}$ is not trivial, then it fixes 2 points $x_1,x_2$ on $C$, maybe conjugate over $k$.  A surface with $e = 1$ is not $G$-minimal. If $e \ge  2$, the surface $X$ has a unique section $S$ with self-intersection $-e$. The fibres $F_1,F_2$ over $x_1,x_2$ are $G$-invariant and each contains  a fixed point on $S$. Since $\ell\ne \cha~k$, the group $G$ must have another fixed point $y_i\in F_i, i = 1,2$. The closed subscheme $\{y_1,y_2\}$ is defined over $k$. After we perform a 
 $G$-invariant elementary transformation at these points, we obtain a surface isomorphic to $\bfF_{e-2}$. Continuing in this way, we arrive at either $\bfF_0 \cong \bbP_k^1\times \bbP_k^1$, or $\bfF_1$, the latter being  not $G$-minimal.  
 \end{proof}

\smallskip
\begin{proposition}\label{prop2} Assume that $\sigma$ of prime order $\ell\ge 7$ acts  minimally on a $k$-rational Del Pezzo surface $X$ of degree $d$. Then one of the following cases occurs:
\begin{itemize}
\item[(i)] $d = 6, t_\ell = 6$;
\item[(ii)] $d = 8, t_\ell = 4$;
\item[(iii)] $d = 9, t_\ell \le 3$. 
\end{itemize} In all cases $X$ has a structure of a toric surface and   $\sigma$ belongs to $T(k)$, where $T$ is an open subset of $X$ isomorphic to a $k$-torus. 
\end{proposition}

\begin{proof} Assume first that  $d\le 5$.  In this case $\Aut(\bar{X})$ is isomorphic to a  subgroup of the Weyl group of a root system of type $A_4, D_{12}, E_6,E_7,$ or $E_8$ (\cite{DI}, Lemma 5.2). The classification of elements of finite order in Weyl groups shows that $\ell \le 7$, and the equality may occur only  if  $d = 2$ or $d = 1$. 

Assume $d = 2$.  Let $\pi:\bar{X}\to \bbP_{\bar{k}}^2$ be a birational morphism that blows down  seven   disjoint $(-1)$-curves  $C_1,\ldots,C_7$ and  $e$ be the divisor class of their sum.  Obviously $e$ has only one effective representative. Let $\calR = (\bbZ K_{\bar{X}})^\perp$ considered as a quadratic lattice isomorphic to the root lattice of type $E_7$. It is known that the  stabilizer subgroup of $e$ in the Weil group $W(E_7)$ is isomorphic to the permutation group $S_7$ (see \cite{Do}). The index of this subgroup is equal  to 576. Elements of the  $W(E_7)$-orbit of $e$ correspond to sets of seven disjoint $(-1)$-curves on $\bar{X}$. They are paired into 288 pairs, two sets in the same pair  differ by the Geiser involution on $\bar{X}$. Recall that the linear system $|-K_{\bar{X}}|$ defines a degree 2 finite map  $\bar{X}\to  \bbP_{\bar{k}}^2$ ramified over a nonsingular plane quartic. The cover transformation is the Geiser involution. It sends $C_i$ to the curve $C_i'$ such that $C_i+C_i'\in |-K_{\bar{X}}|. $ Since $288 \equiv 1 \mod 7$, a group of order 7 has a fixed point on  the set of such pairs and acts identically on each element in the fixed pair.  Thus we can find a $\sigma$-invariant  class $e$ as above. It is known that  $\calR^{\la \sigma\ra} \cong \bbZ$ (see \cite{DI}, Table 5).  Since $ e\cdot K_{\bar{X}} = -7$, the divisor class $2e+7K_{\bar{X}}$ spans $\calR^{\la \sigma\ra}$. Since the Galois group $\Gamma$ commutes with $\la \sigma \ra$, it either leaves the class $2e+7K_{\bar{X}}$ invariant or contains an element that changes $2e+7K_{\bar{X}}$ to $-2e-7K_{\bar{X}}$. In the first case $e$ is  $\la\sigma\ra\times \Gamma$- invariant contradicting the assumption that $(X,\la \sigma\ra)$ is minimal. In the second case $\Gamma$ contains an element  $\gamma$ that sends $e$ to 
$e' = -e-7K_{\bar{X}}$. This implies that $\gamma(C_1+\cdots+C_7) =  \gamma(C_1)+\cdots+\gamma(C_7)$, where 
$C_i+ \gamma(C_i)\in |-K_{\bar{X}}|$. Thus $\gamma$ acts as the Geiser involution. It sends  the divisor class $C_i-C_{j}$ to $C_j-C_i$. Together with $2e+7K_{\bar{X}}$ the divisor classes of $C_i-C_j$  generate $\calR$, hence $\gamma$ acts as the minus identity on $\calR$ and  $\calR^\Gamma = \{0\}$. This implies that $X$ is $k$-minimal, and because $d\le 5$ and $X$ is $k$-rational, this  contradicts Theorem 4.

Assume $d = 1$. An element of order 7 in $W(E_8)$ is conjugate to an element of order 7 in $W(E_7)$, where the inclusion $W(E_7) \subset W(E_8)$ corresponds to the inclusion of the Dynkin diagrams. This 
implies  that $\sigma$ acting on $\bar{X}$ is a lift of an automorphism $\sigma'$ of order 7 of a Del Pezzo surface $\bar{X'}$ of degree 2 under the blow-up of a $\sigma'$-invariant point. Let $C_8$ be the corresponding $\sigma$-invariant exceptional curve. Let $e' = C_1+\cdots+C_7$ be the pre-image of the $\sigma'$-invariant class in $\Pic(\bar{X'})$ from the previous case $d = 2$ and $e = e'+C_8$. Let $\calR = (K_{\bar{X}})^\perp$. Then 
$\calR^{\la\sigma\ra}$ is generated (over $\bbQ$) by $v = \frac{1}{3}(e-3C_8+5K_{\bar{X}})$ and $w = C_8+ K_{\bar{X}}$. We have $(v,v) = -4, (v,w) = 1, (w,w) = -2$.  The Galois group leaves $\calR^{\la\sigma\ra}$ invariant and acts as an isometry of the quadratic form $-4x^2+2xy-2y^2$. It is easy to compute the orthogonal group of this quadratic form to obtain that an element of finite order in this group is either the identity or the minus the identity. In the first case we obtain that $X$ is not $\la \sigma\ra$-minimal. In the second case, by using the argument from the previous case,  we obtain that $\Gamma$ contains an element $\gamma$ that acts as the Bertini involution. Recall that $|-2K_{\bar{X}}|$ defines a degree 2 map $\bar{X}\to Q$, where $Q$ is quadric one in $\bbP_{\bar{k}}^3$. Its cover transformation is the Bertini involution. It   sends $C_i$ to $C_i'\in |-2K_{\bar{X}}-C_i|$. It acts as the minus identity on $\calR$ and hence $X$ is $k$-minimal. We conclude as in the previous case.

 Finally, it remains to consider the case $d \ge 7$. If $d = 7$, the surface $\bar{X}$ is obtained by blowing-up two points in $\bbP_k^2$ and is obviously is not $G$-minimal (the proper transform of the line joining the two points is $G$-invariant and is defined over $k$). If $d = 8$, the surface is either $\bfF_0$ or $\bfF_1$. The surface $\bfF_1$ is obviously not $G$-minimal.  Let $\bar{X} = \bfF_0$. Since $\ell$ is odd, it preserves the rulings on $\bar{X} \cong (\bbP_{\bar{k}}^1)^2$, i.e. the Galois group $\Gamma_k$ switches the two rulings. Obviously, $\sigma$ has 4 fixed points on $\bar{X}$, the vertices of a quadrangle of lines. Since $\sigma$ is a $k$-automorphism, the set of these points defines  a closed subset of $X$ and the quadrangle descends to a divisor $D$ on $X$ whose complement is a $k$-torus $T$. So we obtain that $X$ is a toric Del Pezzo surface of degree 8.  As we saw in Remark \ref{toric}, this case is realized only if $t_\ell= 4$. 
  
 If $d = 9$, the surface  $\bar{X}\cong \bbP_{\bar{k}}^2$. By definition, $X$ is a Severi-Brauer variety of dimension 2. Since $X(k)\ne \emptyset$ it must be trivial, i.e. $X\cong \bbP_k^2$. Thus $\Aut(X) \cong \PGL_3(k)$.  Applying  Theorem \ref{serre1} , we obtain that this can be realized only if $t_\ell \le 3$. The fixed locus of  $\sigma$ in $\bbP_{\bar{k}}^2$ is either the union of a line and a point, or the set of three distinct points. In the first case the line and the point must be defined over $k$. By joining the point with two $k$-points on the line we obtain a triangle of $\sigma$-invariant lines. Its complement is a $k$-torus on which $\sigma$ acts as a translation. In the second case, the three fixed points are fixed under $\sigma$ because $\ell \ge 5$. They defines a $\sigma$-invariant triangle on $\bbP_k^2$ whose complement is a k-torus. 
  
 The remaining case $d = 6$ was considered in Proposition \ref{dp6} from the previous section.  An automorphism of $\sigma$ of order $\ell \ge 5$ acts on a unique Del Pezzo surface $X$ and its image 
 $\Aut(X) \cong T(k)\rtimes D_{12}$ is contained in $T(k)$. This proves Theorem 1.
 \end{proof}
 
 \begin{remark} It is shown in \cite{DI}, Table 8,  that a Del Pezzo surface $X$ of degree 1 over $\bbC$ does not admit an automorphism of order 7, minimal or not. One can show that this is true over an arbitrary   field $k$ (even of  characteristic $7$). 
 
One can also prove the non-existence of a minimal action of a cyclic group of order 7 on a Del Pezzo surface of degrees $d = 1$ or $2$ by using the ideas of Mori theory. We give the proof only in the case $d = 2$ leaving the case $d = 1$ as an exercise to the reader.  Let $G = \la \sigma\ra$. The group $G\times \Gamma$ acts on the Mori cone $\NE(\bar{X})$. It is a polyhedral cone generated by the divisor classes of $56$  $(-1)$-curves on $\bar{X}$. The intersection of  $\NE(\bar{X})$ with $\Pic(\bar{X})^G \cong \bbZ^2$ has two extremal rays, both $G$-invariant. They are switched by the Geiser involution. The surface $\bar{X}$ does not admit a $G$-invariant structure of a conic bundle. In fact,  since $K_X^2 = 2$,  the number of singular  fibres of such a conic bundle must be equal to 6. So $G$ fixes them and switches the components of some of these fibres. Since $G$ is of odd order this is impossible.  Now, it follows from  the proof of Theorem 5 that each $G$-invariant extremal ray is generated by a $\Gamma$-orbit of a $(-1)$-curve that consists of disjoint curves.  Since the number of $(-1)$-curves is equal to $56$ and is divisible by $7$, such an orbit must consist of 7 curves (if the orbit is a singleton, there are at least  $6$ more invariant $(-1)$-curves contradicting to the fact that  $\rank~\Pic(\bar{X})^G = 2$). Thus the $G$-invariant Mori cone is generated by two extremal rays corresponding to two 7-tuples of disjoint  $(-1)$-curves. Since $X$ is $\Gamma-G$-minimal, the group $\Gamma$ switches the two sets sending each curve $C_i$ from the first set to the curve $C_i'\in |-K_X-C_i|$ from the second set. Therefore it acts on $\Pic(X)$ as the Geiser involution. Thus $X$ is $k$-minimal contradicting its rationality. 
\end{remark}
\section{The case of characteristic $0$}

Assume $\cha~k = 0$. We assume also that 

\begin{itemize}
\item[(*)] $k\cap \bbQ(\zeta_\ell) = \bbQ.$
\end{itemize}
Thus 
$$t_\ell = \ell-1.$$

Assume $\ell\ge 7$. By Proposition \ref{prop1}, $\sigma$ cannot act minimally on a conic bundle. By Proposition \ref{prop2}, $\sigma$ can only act minimally on a Del Pezzo surface $X$ of degree 6 in which case $\ell = 7$. By Proposition \ref{dp6},  $X$ is a unique (up to isomorphism) toric surface $T$ over $k$ split over $E = k(\zeta_7)$.   The Galois group $\Gamma$ acts on $T_E$ via the subgroup $H \subset \Aut(\Sigma)$ isomorphic to the cyclic group of order 6. The action of $H$ on $T_E$ is $\Gamma$-equivariant, and hence admits a  descent to an action of $H$ on $X$. The 7-torsion subgroup $T(k)[7]$  of $T(k)$ is $H$-invariant. Hence $H$ acts on the cyclic group $\la \sigma\ra$ of order 7 by automorphisms. This shows that all non-trivial powers of $\sigma$ are conjugate in $\Cr_2(k)$. 

This proves the following. 

\begin{theorem} Assume (*) is satisfied. Then $\Cr_2(k)$ does not contain elements of prime order $> 7$ and all elements of order 7 are conjugate.
\end{theorem}

It follows from the proof of Proposition \ref{prop2} that an element of order $\ell = 5$ can be realized as an automorphism of a Del Pezzo surface of degree 8 defined over $\bbQ$. One can also show that it can be realized as a minimal automorphism of a Del Pezzo surface of degree 5 over $\bbQ$ that arises from the Cremona transformation
$(x,y,z) \mapsto (xz,x(z-y),z(x-y))$ with 4 fundamental points $(1,0,0), (0,1,0),(0,0,1), (1,1,1)$. The blow-up of these points is a Del Pezzo surface of degree 5 (see \cite{DF}, p. 20). Note that $\PGL_3(\bbQ)$  does not contain elements of order $\ge 5$, so these two transformations are not conjugate to a projective transformation. An explicit example of a Cremona transformation of order 7 over $\bbQ$ was given by N. Elkies in his letter to J. P. Serre on November 10, 2005 (see http://math.harvard.edu/$\sim$elkies/serre.pdf).

\bibliographystyle{amsplain}

\end{document}